\documentclass[11pt]{article}
\usepackage[margin=1in]{geometry}
\usepackage[T1]{fontenc}
\usepackage{lmodern}
\usepackage{amsmath,amssymb,amsthm,mathtools}
\usepackage{needspace}
\usepackage{booktabs,microtype,xcolor,tikz,array,tabularx,longtable}
\usetikzlibrary{arrows.meta,positioning}
\usepackage[colorlinks=true,linkcolor=blue!45!black,citecolor=blue!45!black,urlcolor=blue!45!black]{hyperref}
\hypersetup{pdftitle={Adhesive riffle shuffles: one-parameter models and mixing bounds},pdfauthor={Yunjiang Jiang}}
\newtheorem{theorem}{Theorem}[section]
\newtheorem{lemma}[theorem]{Lemma}
\newtheorem{proposition}[theorem]{Proposition}
\newtheorem{corollary}[theorem]{Corollary}
\theoremstyle{definition}
\newtheorem{definition}[theorem]{Definition}
\newtheorem{remark}[theorem]{Remark}
\newcommand{\E}{\mathbb E}
\newcommand{\Prob}{\mathbb P}
\newcommand{\Sn}{\mathfrak S_n}

\newcommand{\mix}{\mathrm{mix}}
\newcommand{\id}{\mathrm{id}}
\newcommand{\Var}{\operatorname{Var}}
\newcommand{\Ent}{\operatorname{H}}

\newcommand{\inner}[2]{\langle #1,#2\rangle}
\newcommand{\norm}[1]{\lVert #1\rVert_2}

\newcommand{\RT}{\mathrm{RT}}
\newcommand{\Card}{\mathsf C}
\newcommand{\Block}{\mathsf B}
\newcommand{\Interp}{\mathsf I}
\newcommand{\Markov}{\mathsf M}
\newcommand{\Lam}{\Lambda_{n,\varepsilon}}
\newcommand{\logplus}{\log_+}
\newcolumntype{Y}{>{\raggedright\arraybackslash}X}
\title{Adhesive riffle shuffles:\\
simple one-parameter models, interpolation,\\
and parameter-dependent mixing bounds}
\author{Yunjiang Jiang}
\date{October 18, 2026}

\begin{document}
\maketitle

\begin{abstract}
We compare four simple shuffles with one adhesion parameter $p$: a
card-weighted sticky riffle, a block-weighted riffle, a midpoint-interpolation
shuffle, and a correlated-label inverse riffle. Independent surviving
contacts create contiguous blocks; no force or energy constants are fitted.
We give polynomial upper bounds with explicit dependence on $p$, common
adjacency-survival lower bounds, and a sharper information lower bound for
interpolation. For fixed deck size $n\ge4$, the mixing time diverges as
$\Theta((1-p)^{-1})$ for the three cut-first models, but as
$\Theta((1-p)^{-2})$ for the correlated-label model. The latter distinction
comes from whether one or two additional contact breaks are needed to
escape cyclic rotations. For even $n$, unglued interpolation retains the
lower bound $\Omega(\sqrt n\log n)$ and upper bound
$O(n^{5/2}\log n)$. We distinguish proved bounds from conjectured scales,
recover classical riffle cutoff under sufficiently small perturbations,
and explain precisely where established results on clumpy shuffles apply.
The fixed-parameter asymptotic orders of the new cut-first families and
their possible cutoff remain unresolved.
\end{abstract}

\noindent\textbf{Status and scope.}
This is an AI-assisted working manuscript, not a refereed publication.
The bounds proved here are the bounds established in this draft, not claims
of optimality in $n$ or of best-known results for every related shuffle.
Standard spectral inputs and external mixing theorems are identified.
Finite computations are diagnostics, not proofs of asymptotic statements.
The historical survey is selective, checked through October 18, 2026;
there is no claim of priority, experimental calibration, or exhaustive
literature coverage. The parameter $p$ is fixed; latent contact and label
randomness is refreshed every round. Persistent electrical charge is
outside the Markov models studied here.

\tableofcontents

\section{Introduction and historical context}

\subsection{From a physical shuffle to a random walk}

The question ``how many shuffles suffice?'' is not well posed until the
shuffle law, the notion of randomness, and the distance from a target law
have been specified. A labelled deck has $n!$ possible orders. A
label-independent shuffling rule is a random walk on $\Sn$, and uniform
measure is stationary because each individual move is a permutation.
Stationarity alone says neither that the walk reaches every order nor how
quickly it forgets its initial state. Modern treatments organize these
questions through mixing times, coupling, stopping times, entropy, and
spectral analysis; see Aldous--Diaconis \cite{AD86} and Levin--Peres
\cite{LP17}. Diaconis--Fulman \cite{DF23} give a broad account of the
mathematics of card shuffling.

Two benchmarks illustrate why the mechanism matters. The random-transposition
walk chooses positions throughout the deck, while adjacent-transposition
walks move cards locally. Diaconis--Shahshahani \cite{DS81} established sharp
mixing for random transpositions using representation theory. For local
shuffling, Wilson \cite{Wilson04} developed sharp-scale bounds and Lacoin
\cite{Lacoin16} proved cutoff for the adjacent-transposition shuffle.
The time conventions in these papers differ: one discrete adjacent update
per round is not the same clock as independent rate-one updates at every
edge. We will specify our clocks whenever comparing chains.

\subsection{Riffle shuffling and cutoff}

The classical Gilbert--Shannon--Reeds (GSR) model draws
$K\sim\operatorname{Bin}(n,1/2)$, cuts off the first $K$ cards, and then
chooses uniformly among the $\binom nK$ interleavings preserving both packet
orders. Its equivalent binary-word description gives $n$ independent fair
input bits. The historical model and equivalent constructions are discussed
by Aldous--Diaconis \cite{AD86}, who trace the model to Gilbert and Shannon
in 1955 and Reeds in 1981, and by Diaconis--Fulman \cite{DF23}.
Here ``Reeds'' refers to Jim Reeds, not to Reed--Solomon error-correcting
codes. In the forward construction the next packet is selected by its
\emph{remaining card count}, not by a fixed fair coin.

Bayer--Diaconis \cite{BD92} analyzed the distribution through rising
sequences and established the sharp transition at
$\tfrac32\log_2 n$ shuffles. The familiar slogan about seven shuffles is a
finite-deck summary, not a statement that seven steps produce exact
uniformity. The distinction between gradual convergence and an asymptotically
abrupt transition is the cutoff phenomenon.

Card-shuffling problems helped make cutoff a central question in finite
Markov-chain theory. The random-transposition analysis of
Diaconis--Shahshahani \cite{DS81} and the stopping-time approach of
Aldous--Diaconis \cite{AD86} supplied complementary algebraic and
probabilistic methods. Diaconis's survey \cite{D96} places shuffling
alongside urn models and other walks, including examples with no cutoff.
Polynomial mixing, a spectral gap, or a plausible entropy threshold
alone does not establish a sharp transition.

Lalley \cite{Lalley99} connected riffle permutations to symbolic and
interval dynamical systems and related randomization to entropy.
This viewpoint motivates studying geometric interleaving rules, but
random point configurations in a geometric representation must not be
confused with resetting packets to deterministic equally spaced grids.
The two constructions need not have the same law.

Recent work of Sellke, Shi, and Wang \cite{SSW25} proves cutoff for broad
classes of cut-size distributions. Their interleaving is still uniformly
random conditional on the cut. Thus randomness of the cut alone is not the
assumption that distinguishes their setting from ours: it is the
conditional interleaving law. Their theorem does not supply an upper bound
for the chain studied here.

\subsection{Perfect shuffles, Thorp shuffling, and sparse randomness}

A perfect shuffle interlaces equal packets in a prescribed alternating
pattern. Such permutations have rich algebraic structure but do not, by
themselves, constitute a mixing mechanism. Diaconis, Graham, and Kantor
\cite{DGK83} study the groups generated by perfect shuffles. Thorp's work on
nonrandom shuffling and faro \cite{Thorp73} is an early motivation for
distinguishing realistic or structured shuffling from uniform random
permutations.

In the Thorp shuffle, the deck is cut exactly in half, and each corresponding
pair of cards is ordered by an independent fair coin. Conditional on an
exactly balanced cut, our shuffle is precisely this operation. Morris's
first polynomial-in-$\log n$ bound appeared in work originating at STOC 2005
and published in 2008 \cite{Morris08}. Morris \cite{Morris09} obtained
$O(\log^4 n)$ for arbitrary even $n$ by an entropy method for shuffles
built from random pair exchanges. For $n=2^d$, Morris \cite{Morris13}
subsequently obtained $O(d^3)$.

An October 3, 2026 preprint of Jiang \cite{Jiang26} reports a full-deck
upper bound of order $d^{5/3}$ for the power-of-two Thorp shuffle.
This recent preprint is included for bibliographic currency; its proof is
not audited or used here. In particular, our argument does not depend on
the best available Thorp mixing bound.

\subsection{Comparison methods and the present contribution}

Comparison of Dirichlet forms and spectra is a standard way to transfer
mixing estimates between walks on groups; Diaconis--Saloff-Coste
\cite{DSC93} provide a foundational treatment. Our application has a
nonreversible complication: the proportional shuffle is not a reversible
walk, so controlling only the eigenvalues of an additive symmetrization
would not suffice. We instead control its one-step squared-norm deficit
on invariant spectral subspaces.

A second ingredient comes from a different subject. The detectability
lemma, introduced by Aharonov, Arad, Landau, and Vazirani \cite{AALV08},
bounds products of local projections. We use the arbitrary-order version
of Anshu, Arad, and Vidick \cite{AAV16} to compare parallel swap layers
with the sum of their individual swap projections. The state space here
is classical; only the Hilbert-space inequality is borrowed.

For unglued interpolation, the model-specific observations are an exact
tie count and a two-cut geometric identity. The identity turns the balanced-cut matching and one
conjugate matching into a Hamiltonian cycle. It permits a polynomial
upper bound without waiting for a long uninterrupted run of balanced cuts.
We also obtain an information-concentration lower bound, stronger than a
fixed-error entropy inequality. The substantial gap between the two bounds
is part of the conclusion, not a suppressed technical issue.

For the adhesive extensions, local alternatives supply fair-swap
projections along a path. A separate comparison uses rotations and
one-break insertions when adhesion is strong. The model atlas and bounds
come first; the unglued proof develops the common spectral tools;
Sections~\ref{sec:adhesive-lower}--\ref{sec:small-adhesion} then give the
parameter-dependent arguments. This organization distinguishes baseline
results, uniform bounds, and genuinely sharp fixed-deck parameter exponents.

\subsection{Physics motivates blocks, not an uncalibrated force law}

Surface-charge correlations can contribute to adhesion between insulating
surfaces \cite{PSV13}. Experiments on insulating polymers also link charge
to friction \cite{Burgo13}, while the mechanics of interleaved paper exhibits
strong geometry-dependent friction amplification \cite{Alarcon16}.
These are motivations for allowing coherent groups of cards, not
measurements of a playing-card shuffle. Sticking need not be electrical,
and a net positive or negative charge per card does not describe all
surface-charge effects. We do not insert a thermal Boltzmann factor into
a macroscopic, hand-driven process.

Our coarse-graining keeps just one probability: $p$ is the chance that an
existing neighboring contact survives the release attempt. Broken contacts
permit, but do not force, separation in the output. Humidity, coating,
pressure, and handling would affect an experimentally fitted $p$; they are
not extra mathematical parameters here. Independent contacts, equal card
thickness, and complete refresh between shuffles are simplifying
assumptions, not consequences of electrostatics. The cut-first models
also impose the cut independently of adhesion. Persistent charge would
require an enlarged state space and a new stationarity argument.

\subsection{Clumpy shuffles and the cut-law distinction}

Jonasson and Morris \cite{JM15} prove $O(\log^4 n)$ mixing for a riffle
driven by a stationary, symmetric two-state Markov word, with any fixed
nondegenerate transition probability. Their paper cites earlier work by
Wager on an extremely clumpy regime; we have not independently checked
that thesis and make no novelty claim in this regime. Our correlated-label
family below is exactly their model, up to inversion and parameterization.
Our other three families keep a binomial cut and are different chains.
In particular, their polylogarithmic theorem is not an upper bound for
those cut-first families.

Contiguous random packets also appear in the overhand-shuffle literature.
Pemantle \cite{Pemantle89} proved bounds between orders $n^2$ and
$n^2\log n$; Jonasson \cite{Jonasson06} established the latter as the
correct order for that model. Its packet-reversal mechanism is different
from preserving two packets while interleaving their blocks. Palomares's
2025 REU report \cite{Palomares25} studies the Markovian rule under the
name ``geometric model,'' with a neatness parameter and finite-deck
spectral calculations. That use of ``geometric'' refers to run lengths,
not our deterministic midpoint grids. Neither resemblance transfers a
mixing theorem to the cut-first models.

\section{A model atlas: one adhesion parameter, four rules}
\label{sec:atlas}

Throughout the extension let $p\in[0,1]$ and $q=1-p$. These are not two
parameters. Card identities never affect the rule. All random choices
are fresh at each shuffle. The parameter $p$ is contact survival, not the
probability that two successive output cards come from the same packet.

\subsection{The common cut-first mechanism}

The families $\Card_p,\Block_p,\Interp_p$ first draw
$K\sim\operatorname{Bin}(n,1/2)$ and form packets of $K$ and $n-K$ cards.
Empty cuts are the identity. At each internal packet gap independently,
keep the bond with probability $p$. The connected components are blocks;
each block keeps its internal order. Conditional on $0<K<n$,
\begin{equation}\label{eq:block-count}
 B=2+\operatorname{Bin}(n-2,q),\qquad \E B=2+(n-2)q.
\end{equation}
Before boundary truncation, block lengths have law
$\Prob(L=\ell)=q p^{\ell-1}$ and mean $1/q$ for $p<1$.

\begin{definition}[Card-weighted release: $\Card_p$]\label{def:card}
With $A,B$ cards left in the packets, choose the next left block with
probability $A/(A+B)$ and otherwise the next right block. Release the
whole selected block and subtract its length from the packet count.
Equivalently, after each individual card immediately update the counts.
If both packets remain nonempty, test the next bond on the last-used
side: with probability $p$ continue on that side, and with probability
$q$ make a fresh card-weighted GSR choice.
Thus, with $A_t,B_t$ measured \emph{after} $t$ drops and $H_t$ the last
side, the next-left probability is
\begin{equation}\label{eq:updated-counts}
 p\mathbf1_{\{H_t=L\}}+q\frac{A_t}{A_t+B_t}
 \quad (A_tB_t>0).
\end{equation}
An empty packet cannot supply another card; the first drop has no
previous side and uses $K/n$. A failed bond can still be followed by a
fresh choice of the same side.
\end{definition}

\begin{definition}[Block-weighted release: $\Block_p$]\label{def:block}
With $r,s$ blocks remaining, choose the next left block with probability
$r/(r+s)$. Conditional on the packet partitions, all order-preserving
block interleavings are equally likely. Blocks of unequal lengths have
equal selection weight. This is deliberately \emph{not} card-weighted
GSR. Nor is its cut binomial in the number of blocks: the initial cut was
made in cards.
\end{definition}

\begin{definition}[Midpoint release: $\Interp_p$]\label{def:sticky-linear}
A block at packet positions $a,\ldots,b$ in a packet of $m$ cards has key
$(a+b-1)/(2m)$. Merge blocks in increasing key order. Independent fair
coins resolve equal keys, and each block is emitted intact. The keys use
card positions, not block ranks. This is a release schedule approximation,
not an independently calibrated mechanical clock.
\end{definition}

\begin{figure}[ht]
\centering
\begin{tikzpicture}[node distance=5mm and 8mm,
 box/.style={draw,rounded corners,align=center,font=\small,inner sep=5pt},
 >={Latex[length=2mm]}]
\node[box] (cut) {Binomial cut\\in \textbf{cards}};
\node[box,right=of cut] (bonds) {Independent bonds\\survive with probability $p$};
\node[box,right=of bonds,yshift=16mm] (card) {$\Card_p$: remaining\\\textbf{card} counts};
\node[box,right=of bonds] (block) {$\Block_p$: remaining\\\textbf{block} counts};
\node[box,right=of bonds,yshift=-16mm] (interp) {$\Interp_p$: block\\midpoint keys};
\draw[->] (cut)--(bonds);
\draw[->] (bonds)--(card);
\draw[->] (bonds)--(block);
\draw[->] (bonds)--(interp);
\end{tikzpicture}
\caption{The three cut-first families share precisely the same cut and
bond law. Only the block interleaver changes.}
\label{fig:atlas}
\end{figure}
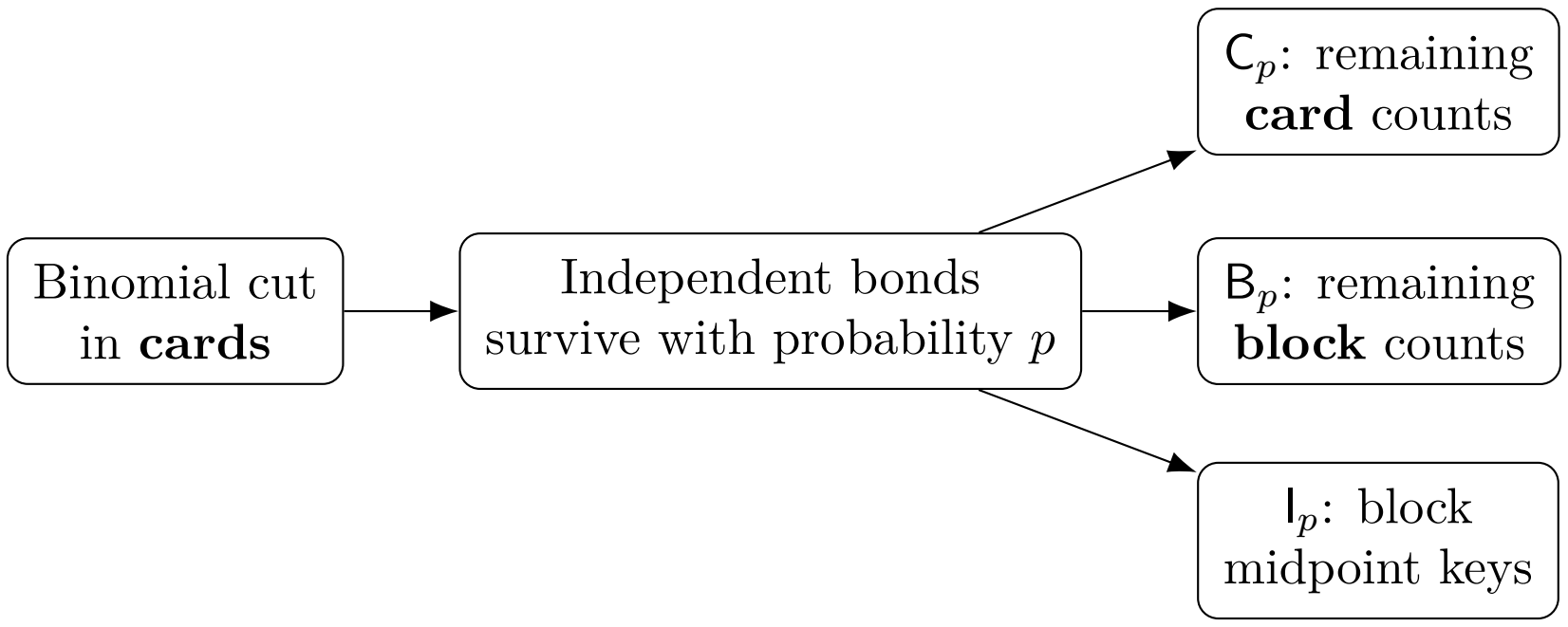

For example, consider packet blocks $[1,2,3,4]$ and $[5],[6],[7],[8]$.
The first-left probability is $1/2$ for $\Card_p$ and $1/5$ for $\Block_p$.
The rule $\Interp_p$ instead gives $5,6,1,2,3,4,7,8$ deterministically.
This example is conditional on a partition; it introduces no new parameter.

\subsection{The correlated-label benchmark has a different cut}

\begin{definition}[Markov inverse riffle: $\Markov_p$]\label{def:markov}
On the uncut deck keep each neighboring bond independently with probability
$p$, give each resulting block an independent fair binary label, and stably
place label $0$ before label $1$. Equivalently, start with a fair label and,
at each next position, repeat the label with probability $(1+p)/2$.
\end{definition}

The label-switch indicators in this construction are independent
Bernoulli$(q/2)$. In the corresponding forward shuffle the number of zeros
is the cut size. Its variance is
\begin{equation}\label{eq:markov-cut}
 \Var(K)=\frac n4+\frac12\sum_{d=1}^{n-1}(n-d)p^d,
\end{equation}
because label correlations at distance $d$ equal $p^d$. Except at $p=0$,
this cut is not binomial. Reversing iid permutation increments preserves
the distance to uniform, so forward and inverse versions have identical
mixing times. This is the connection with \cite{JM15}.

\begin{table}[ht]
\centering\small
\begin{tabularx}{\textwidth}{@{}lYYY@{}}
\toprule
Family & What is randomized? & At $p=0$ & At $p=1$\\
\midrule
$\Card_p$ & Card-weighted block choices & GSR & Rotations only\\
$\Block_p$ & Uniform block interleaving & GSR & Rotations only\\
$\Interp_p$ & Bonds and midpoint tie coins & Original interpolation & Rotations only\\
$\Markov_p$ & Correlated binary labels & Inverse GSR & Identity only\\
\bottomrule
\end{tabularx}
\caption{All four families have at most one adjustable parameter. The
infinitesimal continuous jitter considered later is an equivalent
implementation of $\Interp_0$, not a fifth tunable model.}
\label{tab:models}
\end{table}

Every family is a walk with iid, label-independent permutation increments,
so uniform measure is stationary. The upper bounds below prove convergence
for $p<1$ on even decks, and for $0<p<1$ on odd decks as well. At $p=0$,
the GSR families mix on either parity; our upper arguments leave odd-deck
interpolation untreated. At $p=1$, the rotation subgroup is proper when
$n\ge3$, so stationarity does not imply convergence from a fixed deck.

\subsection{Reading the bounds}
\label{sec:overview}

Use total variation from a specified starting order, uniform target
measure, and one complete shuffle per time unit. Write
\[
 \Lam=\log n+\log(1/\varepsilon),\qquad
 L(n,p)=\frac{1+\logplus(nq)}q,\qquad
 \logplus x=\max\{0,\log x\}.
\]
At $p=1$ mixing on $\Sn$ fails for $n\ge3$. The following table concerns
$p<1$. Constants in $O$ and $\Omega$ are absolute unless a subscript says
otherwise. Lower bounds concern $\varepsilon=1/4$; upper bounds allow
$0<\varepsilon\le1/2$. The polynomial cut-first upper bounds assume
even $n\ge4$; the common rotation bound below also covers odd $n\ge4$
when $0<p<1$.

\begin{table}[ht]
\centering\small
\begin{tabularx}{\textwidth}{@{}lYY@{}}
\toprule
Family & Lower bound established here & Parameter-explicit polynomial upper bound\\
\midrule
$\Card_p$ & $\Omega(L(n,p))$ & $O(n^{7/2}q^{-6}\Lam)$\\[2pt]
$\Block_p$ & $\Omega(L(n,p))$ & $O(n^3(1+npq)q^{-6}\Lam)$\\[2pt]
$\Interp_p$, $p>0$ & $\Omega(\max\{L(n,p),n\log n/\bar h_{n,p}\})$ &
$O\!\left(\dfrac{n^{5/2}\Lam}{q^4\min\{p^2,q^2\}}\right)$\\[5pt]
$\Interp_0$ & $\Omega(\sqrt n\log n)$ & $O(n^{5/2}\Lam)$\\[2pt]
$\Markov_p$ & $\Omega(\max\{L(n,p),(nq^2)^{-1}\})$ & $O(n^2q^{-2}\Lam)$\\
\bottomrule
\end{tabularx}
\caption{Proved bounds, not proposed sharp scales. The entropy quantity
$\bar h_{n,p}$ is explicit in \eqref{eq:sticky-entropy-upper}. The Markov
bound holds for all $n\ge4$; for fixed $p<1$ the stronger external bound
$O_p(\log^4 n)$ is available from \cite{JM15}.}
\label{tab:bounds}
\end{table}

Three refinements are important when interpreting Table~\ref{tab:bounds}.
First, for all cut-first families and $0<p<1$ we prove
\begin{equation}\label{eq:rot-summary}
 t_{\mix}(\varepsilon)\le
 C n^3\left(\frac{n}{p^{n-2}}+
              \frac{1}{q p^{n-3}}\right)\Lam.
\end{equation}
This is particularly useful when $nq$ is bounded, and must not be used
as a fixed-$p$ polynomial bound. Second, interpolation also satisfies
\begin{equation}\label{eq:singleton-summary}
 t_{\mix}^{\Interp_p}(\varepsilon)
 \le C n^{5/2}q^{-(n-2)}\Lam
 \quad\text{for even }n.
\end{equation}
Take the minimum of all applicable upper bounds, not their sum. Third,
at $p=0$ the classical GSR families have the much sharper
$\tfrac32\log_2 n+O_\varepsilon(1)$ bound and cutoff \cite{BD92}.
The constants in the table are deliberately conservative.

\begin{theorem}[Sharp adhesion exponents at fixed deck size]
\label{thm:parameter-exponents}
Fix $n\ge4$ and $0<\varepsilon\le1/4$. As $p\uparrow1$,
\begin{align*}
 t_{\mix}^{\mathsf F_p}(\varepsilon)
   &=\Theta_{n,\varepsilon}((1-p)^{-1}),
   &&\mathsf F\in\{\Card,\Block,\Interp\},\\
 t_{\mix}^{\Markov_p}(\varepsilon)
   &=\Theta_{n,\varepsilon}((1-p)^{-2}).
\end{align*}
The $n$-dependent constants are not asserted to be sharp. These are
fixed-$n$ parameter asymptotics, not a cutoff assertion.
\end{theorem}

The proof is given after the uniform upper and lower estimates. In the
joint regime $nq\le c$, \eqref{eq:rot-summary} improves the cut-first
upper bound to $O_c(n^3q^{-1}\Lam)$ for sufficiently large $n$.
Its lower counterpart is $\Omega(q^{-1})$. For the Markov model the
corresponding bounds are $\Omega((nq^2)^{-1})$ and
$O(n^2q^{-2}\Lam)$. Large gaps in $n$ remain visible.

\section{The unglued interpolation model and its bounds}

\subsection{Midpoint interpolation}

Number the positions $1,\ldots,n$. Given a cut $k$ with $0<k<n$, write the
two packets as $A_1,\ldots,A_k$ and $B_1,\ldots,B_{n-k}$, in their current
orders. Assign coordinates
\begin{equation}\label{eq:coordinates}
 x_i=\frac{i-1/2}{k},\qquad y_j=\frac{j-1/2}{n-k}.
\end{equation}
Merge by increasing coordinate. Whenever $x_i=y_j$, use an independent
fair coin to decide whether $A_i$ or $B_j$ comes first. Cuts $0$ and $n$
are the identity. Coordinates and random choices are reset at every round.

\begin{definition}
The \emph{binomial-cut proportional shuffle} is the preceding merge with
independent cuts $K\sim\operatorname{Bin}(n,1/2)$ at successive rounds.
The name is descriptive, not a claim that this is established terminology.
\end{definition}

There are at most two cards at any coordinate, and distinct ties involve
disjoint cards. The packet orders are always preserved. The midpoint
convention matters: different endpoint or rounding conventions can yield
different chains and different invariants.

\begin{proposition}[Equivalence with infinitesimal jitter]\label{prop:jitter}
Fix $0<\epsilon_n<n^{-2}$. Add independent
$\operatorname{Uniform}(-\epsilon_n,\epsilon_n)$ perturbations to all
coordinates in \eqref{eq:coordinates} before merging. The resulting shuffle
has exactly the law of independent fair tie-breaking.
\end{proposition}
\begin{proof}
A nonzero difference of cross-packet coordinates is an integer multiple
of $1/(2k(n-k))$, hence has magnitude at least $2/n^2$.
The difference between two perturbations is smaller than $2/n^2$, so no
strict cross-packet comparison changes. Within a packet the spacing is
at least $1/n$, which likewise cannot be reversed. At a tie, exchangeability
of the two continuous perturbations gives either order with probability
$1/2$. Disjoint tied pairs use disjoint noises, giving independence.
\end{proof}

\begin{remark}
Independent two-point noises $\pm\epsilon_n$ leave a tie unresolved with
probability $1/2$. Stable left-first sorting would then favor the left
card with probability $3/4$. One can instead give each tied pair opposite
signs, with one fair coin choosing the signs. That is again exactly our
model. Noise large enough to change non-tied comparisons is a different
model; none of the main results silently assumes such additional noise.
\end{remark}

\subsection{Distances and conventions}

Let $\pi_n$ be uniform measure on $\Sn$ and let $X_t$ start from the
identity deck. Set
\[
 d_n(t)=\frac12\sum_{\sigma\in\Sn}
       |\Prob(X_t=\sigma)-\pi_n(\sigma)|,
 \qquad
 t_{\mix}^{(n)}(\varepsilon)=\min\{t:d_n(t)\le\varepsilon\}.
\]
Every deterministic starting order has the same distance by translation.
Entropies use base-two logarithms; unmarked $\log$ denotes the natural
logarithm. The unglued interpolation results in the next five sections
concern even $n\to\infty$; the extension states its size restrictions
separately.

\Needspace{20\baselineskip}
\begin{theorem}[Main bounds]\label{thm:main}
Let $n$ be even and let $T_n(k)$ be the number of ties at cut $k$. Define
\[
 h_n=\Ent(K)+\E T_n(K),\qquad
 a_n=\frac{\log_2(n!)}{h_n}.
\]
Then:
\begin{enumerate}
\item
$h_n\sim\sqrt{n/(2\pi)}$ and
$a_n\sim\sqrt{2\pi n}\log_2 n$.
\item
For every fixed $\delta\in(0,1)$,
\[
 d_n\bigl(\lfloor(1-\delta)a_n\rfloor\bigr)\longrightarrow1.
\]
\item
For $n\ge4$ and $0<\varepsilon\le1/2$, a universal constant $C$ satisfies
\[
 t_{\mix}^{(n)}(\varepsilon)
 \le C n^{5/2}\bigl(\log n+\log(1/\varepsilon)\bigr).
\]
\end{enumerate}
The constants in the upper bound are not optimized.
\end{theorem}

The third assertion follows from a more general statement. Write $P_k$
for the specified fair-tie merge conditional on cut $k$.

\begin{theorem}[Two-cut criterion]\label{thm:two-cut}
Let $n=2m\ge4$ and let $P=\sum_{k=0}^n q_kP_k$ for an arbitrary cut law
with $\beta=\min(q_m,q_{m+1})>0$. Then
\begin{equation}\label{eq:general-upper}
 t_{\mix}(\varepsilon)
 \le C\frac{n^2}{\beta}
       \bigl(\log n+\log(1/\varepsilon)\bigr),
 \qquad 0<\varepsilon\le1/2.
\end{equation}
No assumption on the remaining cut probabilities is needed.
\end{theorem}

For the binomial law,
$q_m=2^{-n}\binom n m\asymp n^{-1/2}$ and
$q_{m+1}=q_m m/(m+1)$, which gives Theorem~\ref{thm:main}(3).
The case $n=2$ is a two-state chain and is elementary.

The lower bound already confirms a qualitative difference from the
benchmarks: this even-deck chain mixes asymptotically more slowly than
both GSR riffle shuffling and the Thorp shuffle. Indeed, its
$\Omega(\sqrt n\log n)$ lower bound exceeds both the GSR logarithmic
scale \cite{BD92} and the general-even Thorp upper bound
$O(\log^4 n)$ \cite{Morris09}. This comparison is asymptotic, not a
numerical prediction for a 52-card deck.

\section{Tie arithmetic and the entropy of one round}

\begin{lemma}[Exact tie count]\label{lem:ties}
For $0<k<n$, put $g=\gcd(k,n-k)$, $a=k/g$, and $b=(n-k)/g$. Then
\begin{equation}\label{eq:tie-count}
 T_n(k)=
 \begin{cases}
  g,& a\text{ and }b\text{ are both odd},\\
  0,&\text{otherwise}.
 \end{cases}
\end{equation}
Also $T_n(0)=T_n(n)=0$.
\end{lemma}
\begin{proof}
The tie equation is $b(2i-1)=a(2j-1)$. Since $a$ and $b$ are coprime,
both must be odd for a solution. If they are odd, the solutions are
$2i-1=a\ell$ and $2j-1=b\ell$, with
$\ell=1,3,\ldots,2g-1$. There are exactly $g$ solutions.
\end{proof}

For $n=2m$, $T_n(m)=m$. A nonzero tie count requires $g\mid m$, so
\begin{equation}\label{eq:gcd-bound}
 T_n(k)\le\gcd(k,m),\qquad
 T_n(k)\le |k-m|\quad(k\ne m).
\end{equation}
The latter follows because a positive $g$ divides both $k$ and $m$.
If $n$ is odd there are no ties at all. If $n=2^d$, every interior cut
has $T_n(k)=2^{v_2(k)}$.

\begin{lemma}[The balanced cut dominates]\label{lem:entropy}
For $n=2m$ and $p_n=\Prob(K=m)$,
\begin{equation}\label{eq:tie-expectation}
 0\le \E T_n(K)-mp_n
 \le \sum_{r\mid m}\frac{\varphi(r)}r
 \le \tau(m),
\end{equation}
where $\varphi$ is Euler's totient and $\tau$ counts positive divisors.
Consequently,
\[
 h_n=mp_n+O(\log n+\tau(m))\sim\sqrt{\frac{n}{2\pi}}.
\]
\end{lemma}
\begin{proof}
Write $p_j=\Prob(K=m+j)$. The positive-index sequence is decreasing and
$p_j=p_{-j}$. For every divisor $r$ of $m$, grouping $r$ consecutive
positive indices gives
\[
 \Prob(K\ne m,\ r\mid K)
 =2\sum_{\ell=1}^{m/r}p_{\ell r}
 \le \frac{2}{r}\sum_{j=1}^m p_j\le\frac1r.
\]
Use \eqref{eq:gcd-bound} and
$\gcd(k,m)=\sum_{r\mid m,\ r\mid k}\varphi(r)$ to obtain
\eqref{eq:tie-expectation}. Also $\Ent(K)\le\log_2(n+1)$.
The elementary divisor bound $\tau(m)=m^{o(1)}$, or just
$\tau(m)=O(m^{1/4})$, makes the error $o(\sqrt n)$.
For completeness, the latter estimate follows by writing
$m=\prod p^{a_p}$: for $p\ge16$,
$a_p+1\le2^{a_p}\le p^{a_p/4}$, while the finitely many smaller primes
contribute a bounded factor to
$\prod(a_p+1)/p^{a_p/4}$.
Finally Stirling's formula gives
$p_n=2^{-2m}\binom{2m}{m}\sim(\pi m)^{-1/2}$.
\end{proof}

\begin{remark}[Input entropy is not automatically deck entropy]
The effective seed for one step consists of the cut and the tie coins.
Its entropy is exactly $h_n$. Different seeds may produce the same
permutation, so the entropy of the one-step permutation is only bounded
above by $h_n$. After $t$ rounds the deck entropy is at most $th_n$.
There is no assertion that all input information is retained efficiently.
\end{remark}

The exactly balanced cut occurs with probability of order $n^{-1/2}$ and
uses $n/2$ coins. These relatively rare, high-entropy rounds dominate the
leading entropy budget. Thus it would be incorrect to count only
$n+1$ possible moves per round, as one could for deterministic tie rules.

\section{A strong information lower bound}

An entropy inequality alone gives a fixed-error lower bound. Concentration
of the information of the entire input history gives distance tending to
one, which is the stronger conclusion needed for the lower side of a
possible cutoff.

\begin{lemma}[Information moments]\label{lem:information}
Let $q_k=\Prob(K=k)$ and define
\[
 Z_n=-\log_2 q_K+T_n(K).
\]
Then $\E Z_n=h_n$ and $\Var(Z_n)=O(n^{3/2})$ for even $n$.
\end{lemma}
\begin{proof}
Conditional on $K=k$, each tie-coin string has probability $2^{-T_n(k)}$,
which proves the mean assertion. For a variable on $M$ possible values,
its surprisal $A=-\log_2\Prob(K)$ satisfies
$\Prob(A\ge u)\le M2^{-u}$. Integrating this tail with $M=n+1$ gives
$\E A^2=O(\log^2 n)$.
By \eqref{eq:gcd-bound},
\[
 \E T_n(K)^2\le m^2p_n+\E(K-m)^2
 =m^2p_n+n/4=O(n^{3/2}).
\]
Now use $(A+T)^2\le2A^2+2T^2$.
\end{proof}

\begin{proof}[Proof of Theorem~\ref{thm:main}(1)--(2)]
Part (1) follows from Lemma~\ref{lem:entropy} and Stirling's formula for
$n!$. Put $L_n=\log_2(n!)$ and
$t=\lfloor(1-\delta)L_n/h_n\rfloor$.
The information of a full history of cuts and tie coins is
$I_t=\sum_{s=1}^t Z_{n,s}$, a sum of independent copies. Hence
\[
 \E I_t\le(1-\delta)L_n,
 \qquad \Var(I_t)=O(n^2\log n).
\]
Chebyshev's inequality implies
\[
 \Prob\bigl(I_t\le(1-\delta/2)L_n\bigr)
 \ge1-O\bigl((\delta^2\log n)^{-1}\bigr).
\]
Every history in this event has probability at least
$2^{-(1-\delta/2)L_n}$, so there are at most
$2^{(1-\delta/2)L_n}$ such histories. Their resulting decks form a set
$F_n$ with the same cardinality bound. Therefore
\[
 \Prob(X_t\in F_n)\to1,
 \qquad \pi_n(F_n)\le(n!)^{-\delta/2}\to0,
\]
and total variation tends to one.
\end{proof}

\begin{proposition}[A finite-deck entropy bound]\label{prop:finite}
If $d_n(t)\le\varepsilon\le1/2$, then
\[
 t\ge
 \frac{(1-\varepsilon)\log_2(n!)-H_2(\varepsilon)}{h_n},
\]
where $H_2(u)=-u\log_2u-(1-u)\log_2(1-u)$.
\end{proposition}
\begin{proof}
Couple $X_t$ to a uniform permutation with disagreement probability
$d_n(t)$. Conditioning on the disagreement indicator bounds the missing
entropy by $H_2(\varepsilon)+\varepsilon\log_2(n!)$. Combine this with
$\Ent(X_t)\le th_n$.
\end{proof}

\section{Two cuts produce a cycle of random-swap locations}
\label{sec:geometry}

Let $n=2m\ge4$. A permutation $d_k$ will mean the output list of the
midpoint merge at cut $k$, with all ties resolved left before right.
Thus $d_k$ maps an output position to its input position. Permutations
compose as functions: $(uv)(i)=u(v(i))$.

Let $M$ be the balanced input matching
\[
 M=\{\{i,m+i\}:1\le i\le m\}.
\]
Write $H_M$ for the group generated by these disjoint transpositions,
put $b=d_m$, $c=d_{m+1}$, and set $R=cb^{-1}$.

\begin{lemma}[Two-cut geometry]\label{lem:cycle}
The following hold.
\begin{enumerate}
\item The conditional permutation at cut $m$ is $hb$, with $h$ uniform
on $H_M$.
\item At cut $m+1$ there is either no tie or one tie, and any tie is an
edge of $M$.
\item $M$ and $M'=RMR^{-1}$ are disjoint perfect matchings whose union is
a Hamiltonian cycle on the $n$ input positions.
\end{enumerate}
\end{lemma}
\begin{proof}
The first assertion is the definition of the balanced merge.
At cut $m+1$, the packet sizes are $m+1$ and $m-1$. If $m$ is odd,
Lemma~\ref{lem:ties} gives no tie. If $m$ is even, there is exactly one,
between input positions $m/2+1$ and $3m/2+1$, an edge of $M$.

For the last assertion, $b^{-1}Mb$ pairs consecutive output positions.
Thus $M'$ consists of consecutive pairs in the list $c$.
Set $r=\lfloor(m-1)/2\rfloor$. Direct comparison of midpoint coordinates
shows that those pairs are
\begin{align*}
 &\{j,m+1+j\}, &&1\le j\le r,\\
 &\{r+1,r+2\},&&\\
 &\{m+1+j,j+2\},&&r+1\le j\le m-1.
\end{align*}
One way to check the list is to note that $B_j$ lies between $A_j$ and
$A_{j+1}$ when $2j<m$, and between $A_{j+1}$ and $A_{j+2}$ when
$2j>m$. At equality the prescribed left-first convention places
$A_{j+1}$ first.

Contract each edge $\{i,m+i\}$ of $M$ to vertex $i$.
The displayed pairs become exactly
$\{1,2\},\{2,3\},\ldots,\{m-1,m\},\{m,1\}$.
For $m=2$ this is the two-edge cycle as a multigraph.
Hence the original union of the two matchings is connected and
two-regular, and is therefore one cycle of length $n$.
\end{proof}

\begin{figure}[ht]
\centering
\begin{tikzpicture}[scale=1.05,every node/.style={font=\small}]
\foreach \index/\labeltext in {0/1,1/5,2/8,3/4,4/7,5/3,6/2,7/6}
  \node[draw,circle,fill=white,inner sep=3pt] (vertex\index)
    at ({90-45*\index}:1.75) {$\labeltext$};
\foreach \first/\second in {0/1,2/3,4/5,6/7}
  \draw[blue!65!black,very thick] (vertex\first)--(vertex\second);
\foreach \first/\second in {1/2,3/4,5/6,7/0}
  \draw[orange!85!black,very thick,dashed] (vertex\first)--(vertex\second);
\node at (4,0.4) {\textcolor{blue!65!black}{solid: $M$}};
\node at (4,-0.2) {\textcolor{orange!85!black}{dashed: $RMR^{-1}$}};
\end{tikzpicture}
\caption{For $n=8$, the balanced cut $4$ and adjacent cut $5$ produce
the alternating cycle $1,5,8,4,7,3,2,6,1$. This is a geometric identity
between matchings, not an assumption that balanced cuts occur consecutively.}
\label{fig:cycle}
\end{figure}
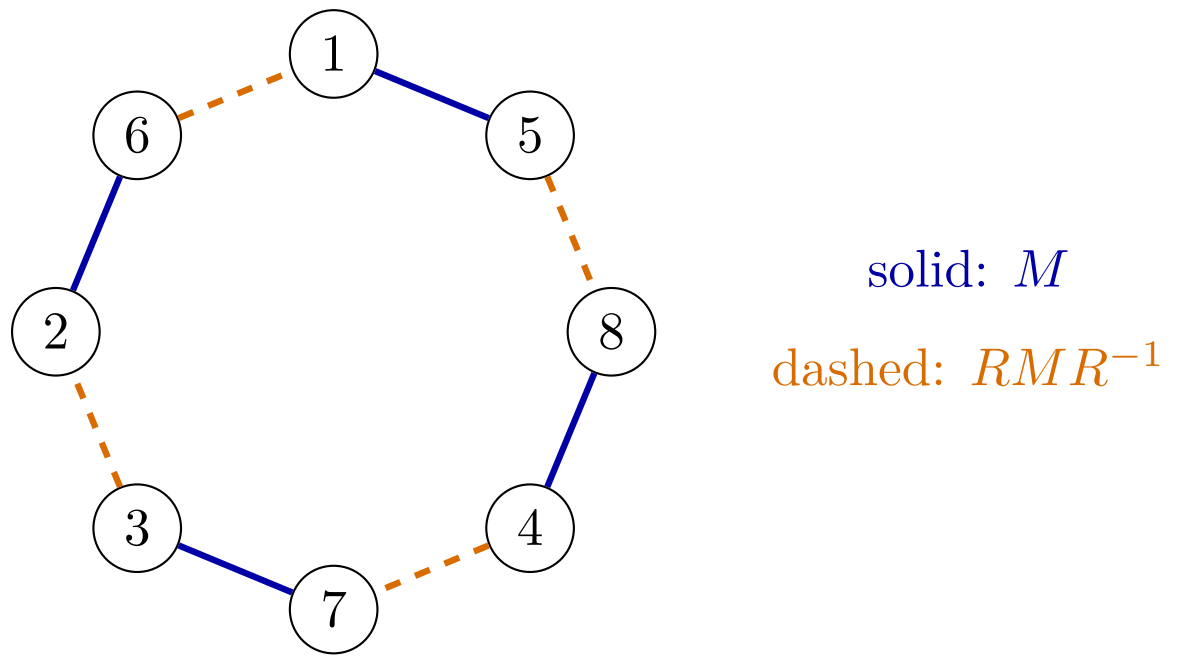

This also proves irreducibility for any cut law charging $m$ and $m+1$:
the group generated by the support contains $H_M$, $b$, and $c$, hence
$R$ and all transpositions in $M\cup M'$. Transpositions along a connected
graph generate $\Sn$. For the binomial chain, the identity has positive
probability, so aperiodicity follows immediately. The quantitative
argument below also handles the two-cut law without an explicit identity
atom.

\section{A polynomial upper bound}

\subsection{Operators and the nonreversible deficit}

Work in $\ell^2(\pi_n)$ and let
$U_gf(\sigma)=f(\sigma g)$. Then $U_gU_h=U_{gh}$ and $U_g$ is unitary.
For a matching $N$, define its averaging projection
\[
 C_N=\prod_{e\in N}\frac{I+U_{\tau_e}}2.
\]
The factors commute; $C_N$ is the orthogonal projection onto invariance
under independently swapping the matched pairs.

Let $P$ be the transition operator in Theorem~\ref{thm:two-cut} and
$Q=P^*$ its adjoint. Put $A=C_M$ and $B=C_{M'}=U_RAU_R^{-1}$.
Conditional on the two distinguished cuts,
\[
 Q_0=U_{b^{-1}}A,\qquad Q_1=U_{c^{-1}}E,
\]
where $E=I$ or $E=(I+U_\tau)/2$ for a transposition $\tau\in M$.
In particular $EA=A$, and
$\norm{(I-E)f}\le\norm{(I-A)f}$.

Define the nonnegative norm deficit
\[
 \mathcal D(f)=\norm f^2-\norm{Qf}^2.
\]
For any mixture $Q=\sum q_kQ_k$ of contractions, the variance identity gives
\begin{equation}\label{eq:mixture-variance}
 \mathcal D(f)=\sum_kq_k\bigl(\norm f^2-\norm{Q_kf}^2\bigr)
       +\sum_kq_k\norm{Q_kf-Qf}^2.
\end{equation}

\begin{lemma}[Control of two parallel layers]\label{lem:layer-control}
With $S=(A+B)/2$ and $\beta=\min(q_m,q_{m+1})$,
\begin{equation}\label{eq:layer-control}
 \inner f{(I-S)f}\le\frac{16}{\beta}\mathcal D(f).
\end{equation}
\end{lemma}
\begin{proof}
Put $a=\norm{(I-A)f}$ and $c_0=\norm{Q_0f-Q_1f}$.
Equation~\eqref{eq:mixture-variance} implies
\[
 a^2\le\mathcal D(f)/\beta,\qquad
 c_0^2\le2\mathcal D(f)/\beta.
\]
For the second inequality, retain the two distinguished terms of the
weighted variance and minimize their common center; their coefficient
is $q_mq_{m+1}/(q_m+q_{m+1})\ge\beta/2$.

Multiplication by $U_c$ shows
$c_0=\norm{U_RA f-Ef}$. Consequently
\[
 \norm{(U_R-I)f}\le c_0+2a,
 \qquad \norm{(U_R-I)f}^2\le12\mathcal D(f)/\beta.
\]
Since $B=U_RAU_R^{-1}$,
\[
 \norm{(I-B)f}
 =\norm{(I-A)U_R^{-1}f}
 \le a+\norm{(U_R^{-1}-I)f}.
\]
Its square is at most $26\mathcal D(f)/\beta$. Averaging this with the
bound for $a^2$ gives \eqref{eq:layer-control} (in fact $27/2$ suffices).
\end{proof}

\subsection{From a cycle to the central transposition generator}

For any transposition edge $e$, let
$F_e=(I-U_{\tau_e})/2$, an orthogonal projection. The cycle Hamiltonian is
\[
 \mathcal H=\sum_{e\in M\cup M'}F_e.
\]
The terminology ``Hamiltonian'' here means only a sum of positive
self-adjoint operators. Define the central random-transposition generator
\begin{equation}\label{eq:RT}
 L_{\RT}=\frac1{\binom n2}\sum_{1\le i<j\le n}(I-U_{(ij)}).
\end{equation}
It generates a continuous-time walk making one uniformly chosen
transposition per unit time on average. The continuous clock avoids the
parity periodicity of a discrete walk that always transposes.

\begin{lemma}[Path comparison]\label{lem:path}
As quadratic forms,
\begin{equation}\label{eq:path}
 L_{\RT}\preceq8n\mathcal H.
\end{equation}
\end{lemma}
\begin{proof}
Delete one cycle edge to obtain a path. A transposition of vertices at
path distance $\ell$ is a word of length $2\ell-1$ in adjacent
transpositions, each path edge used at most twice. Telescoping and
Cauchy--Schwarz give
\[
 \norm{(I-U_{(ij)})f}^2
 \le4n\sum_{e\in\operatorname{path}(i,j)}
          \norm{(I-U_{\tau_e})f}^2.
\]
Each fixed edge lies in at most $n^2/4$ such paths. Average over the
$\binom n2$ pairs, using
$\inner f{L_{\RT}f}=(2\binom n2)^{-1}\sum_{i<j}
\norm{(I-U_{(ij)})f}^2$ and
$\inner f{F_ef}=\norm{(I-U_{\tau_e})f}^2/4$.
The resulting coefficient is at most $8n$.
This is an elementary instance of the path comparison method
\cite{DSC93}.
\end{proof}

The next input prevents losing a further factor of $n$ when replacing
the individual edge projections by simultaneous matching updates.

\begin{lemma}[Projection-product input]\label{lem:DL}
Suppose orthogonal projections $F_e$ each fail to commute with at most
two others. Let $D$ be the product of their complements, each applied
once in any order; this includes any grouping into commuting layers.
On a common invariant subspace $V$ where
$\sum_eF_e\succeq\gamma I$, one has
\[
 \norm{D|_V}^{\,2}\le(1+\gamma/4)^{-1}.
\]
\end{lemma}

This is the detectability-lemma estimate of Anshu--Arad--Vidick
\cite[Lemma 2]{AAV16}, applied with noncommutation degree two. More precisely,
their energy form of the lemma applies on the invariant subspace even
if the restricted sum has no zero eigenvalue. No tensor-product
assumption on $V$ is needed. In our cycle, disjoint transpositions
commute and each edge has just two neighboring edges, so the hypothesis
holds.

\begin{lemma}[Central spectral comparison]\label{lem:central}
For every $f\in\ell^2(\pi_n)$,
\begin{equation}\label{eq:central}
 \norm f^2-\norm{Qf}^2
 \ge \frac{\beta}{16384n}\inner f{L_{\RT}f}.
\end{equation}
Moreover the same inequality holds on every eigenspace of $L_{\RT}$,
and these eigenspaces are invariant under $Q$.
\end{lemma}
\begin{proof}
The sum of all transpositions is central in the group algebra. Thus every
$U_g$, and hence every operator under consideration, preserves each
eigenspace $V_\lambda$ of $L_{\RT}$, whose eigenvalues lie in $[0,2]$.
By Lemma~\ref{lem:path},
$\mathcal H|_{V_\lambda}\succeq\lambda/(8n)\,I$.
Lemma~\ref{lem:DL} implies
\[
 \norm{BAf}^2\le\frac{\norm f^2}{1+\lambda/(32n)},
 \qquad f\in V_\lambda.
\]
In particular,
$\norm f^2-\norm{BAf}^2\ge\lambda\norm f^2/(64n)$.
On the other hand, writing $a=\norm{(I-A)f}$ and
$b_0=\norm{(I-B)f}$,
\begin{align*}
 \norm f^2-\norm{BAf}^2
 &=a^2+\norm{(I-B)Af}^2\\
 &\le3a^2+2b_0^2
 \le6\inner f{(I-S)f}.
\end{align*}
Therefore
$\inner f{(I-S)f}\ge\lambda\norm f^2/(384n)$.
Combining with Lemma~\ref{lem:layer-control} gives
$\mathcal D(f)\ge\beta\lambda\norm f^2/(6144n)$;
the stated constant is weaker. Orthogonal decomposition into the
$V_\lambda$ gives the result for arbitrary $f$.
\end{proof}

\subsection{Why the full spectrum gives the improved bound}

Using only a spectral gap and the smallest stationary mass would lose
an unnecessary factor comparable to $n$. Instead we retain all central
spectral subspaces.

\begin{lemma}[Random-transposition input]\label{lem:RT-mixing}
For the continuous-time chain generated by \eqref{eq:RT}, let
$\chi^2_{\RT}(s)$ be the squared $L^2(\pi_n)$ distance of its density
from $1$, starting at the identity. There is a universal constant $C_0$
such that
\[
 \chi^2_{\RT}(s)\le4\varepsilon^2
 \quad\text{if}\quad
 s\ge C_0n\bigl(\log n+\log(1/\varepsilon)\bigr),
 \qquad 0<\varepsilon\le1/2.
\]
\end{lemma}

\begin{proof}
We give a nonsharp spectral-sum bound to make the required estimate
explicit. The classical representation-theoretic formulas for random
transpositions \cite{DS81} index the spectrum by partitions $\rho$ of
$n$. If $d_\rho$ is the irreducible dimension and $\rho'$ the conjugate
partition, the eigenvalue of the transposition averaging operator is
\[
 r_\rho=
 \frac{\sum_i\binom{\rho_i}{2}-\sum_i\binom{\rho'_i}{2}}
      {\binom n2},
\]
with multiplicity $d_\rho^2$. Consequently
\[
 \chi^2_{\RT}(s)
 =\sum_{\rho\ne(n)}d_\rho^2e^{-2s(1-r_\rho)}.
\]
Put $j=n-\rho_1$. Since all row lengths are at most $n-j$,
\[
 \sum_i\binom{\rho_i}{2}\le\frac{n(n-j-1)}2,
 \qquad 1-r_\rho\ge\frac{j}{n-1}.
\]
Deleting the first row leaves a partition $\mu$ of $j$. A standard
tableau of shape $\rho$ is determined by the $j$ entries below its first
row and a standard tableau of shape $\mu$; column constraints can only
reduce the count. Thus $d_\rho\le\binom nj d_\mu$ and, using
$\sum_{\mu\vdash j}d_\mu^2=j!$,
\[
 \sum_{\rho:\,n-\rho_1=j}d_\rho^2
 \le\binom nj^2j!\le\frac{n^{2j}}{j!}.
\]
It follows that
\begin{equation}\label{eq:RT-explicit}
 \chi^2_{\RT}(s)
 \le\exp\!\left(n^2e^{-2s/(n-1)}\right)-1.
\end{equation}
For $s\ge(n-1)(\log n+\log(1/\varepsilon))$, the right side is at most
$e^{\varepsilon^2}-1\le4\varepsilon^2$ when
$0<\varepsilon\le1/2$. In particular $C_0=1$ suffices.
\end{proof}

The preceding estimate uses the full spectral sum, not merely its
smallest positive eigenvalue. It is weaker than the sharp classical
random-transposition analysis, but suffices for our comparison.

\begin{proof}[Proof of Theorem~\ref{thm:two-cut}]
Set $\alpha=\beta/(16384n)$. By Lemma~\ref{lem:central},
\[
 \norm{Q|_{V_\lambda}}^2\le1-\alpha\lambda.
\]
Let $f_0=n!\mathbf1_{\{\id\}}$ be the density of the initial point mass.
The squared norm of its projection onto $V_\lambda$ is
$\dim(V_\lambda)$. Indeed, the spectral projector commutes with all
translations, so its diagonal is constant and its trace is its rank.
The density after $t$ steps is $Q^tf_0$. Consequently,
\begin{align*}
 \norm{Q^tf_0-1}^2
 &\le\sum_{\lambda>0}\dim(V_\lambda)(1-\alpha\lambda)^t\\
 &\le\sum_{\lambda>0}\dim(V_\lambda)e^{-\alpha t\lambda}
 =\chi^2_{\RT}(\alpha t/2).
\end{align*}
Use $d_n(t)\le\tfrac12\norm{Q^tf_0-1}$ and
Lemma~\ref{lem:RT-mixing}. More explicitly, the estimates give
\[
 t_{\mix}(\varepsilon)\le
 \left\lceil\frac{32768n(n-1)}{\beta}
       \bigl(\log n+\log(1/\varepsilon)\bigr)\right\rceil,
\]
which proves \eqref{eq:general-upper}. This argument establishes convergence even
without separately imposing aperiodicity on the two-cut law.
\end{proof}

\begin{remark}[What is and is not compared]
The proof does not assert that interspersing arbitrary moves can only
accelerate Thorp shuffling, nor does it thin the trajectory to balanced
cuts. Both would require additional justification. It compares the
actual one-step contraction of the full chain, including all intervening
moves, to a central operator. The use of $Q=P^*$, rather than an
unjustified reversible replacement for $P$, is essential.
\end{remark}

\section{Lower bounds for the adhesive families}
\label{sec:adhesive-lower}

\subsection{A universal information floor}

Each cut-first increment is an ordinary order-preserving two-packet
interleaving, irrespective of its probability. There are at most $2^n$
such increments, counted by their binary output words. The same support
bound holds for their inverses, including $\Markov_p$. Consequently,
\begin{equation}\label{eq:support-lower}
 d_n(t)\ge1-\frac{2^{nt}}{n!},\qquad
 t_{\mix}(\varepsilon)\ge
 \frac{\log_2(n!)+\log_2(1-\varepsilon)}n.
\end{equation}
This gives $\Omega(\log n)$ without identifying seed entropy with the
entropy of the output permutation.

\subsection{Circular adjacency is the correct observable}

Let $J(\sigma)$ count original \emph{directed circular} neighbors
$(1,2),(2,3),\ldots,(n,1)$ that remain neighbors in $\sigma$, including
the bottom-to-top edge. Then $J(\id)=n$ and
\begin{equation}\label{eq:adj-uniform}
 \E_{\pi_n}J=\frac{n}{n-1}.
\end{equation}
Circular rather than linear adjacency matters: a rotation leaves $J$
unchanged.

\begin{proposition}[Finite adhesion lower bounds]\label{prop:adj-lower}
Put $r_n=1-2^{1-n}$. For any cut-first family and $n\ge4$,
\begin{equation}\label{eq:cut-adj-lower}
 d_n(t)\ge1-\frac1{n-1}-\frac{3t r_n(n-2)q}{n}.
\end{equation}
For $\Markov_p$ one has
\begin{equation}\label{eq:markov-adj-lower}
 d_n(t)\ge1-\frac1{n-1}
       -\frac{3t(n-1)(n-2)q^2}{8n}.
\end{equation}
In addition, for all four families, every integer $\ell\ge1$ satisfies
\begin{equation}\label{eq:survival-binomial}
 d_n(t)\ge
 \Prob\{\operatorname{Bin}(n,p^t)\ge2t+\ell\}
       -\frac{n}{(n-1)\ell}.
\end{equation}
\end{proposition}

\begin{proof}
For a nonempty cut, let $R$ be the number of broken internal bonds.
It has law $\operatorname{Bin}(n-2,q)$. If $R=0$, the shuffle is a
rotation. Otherwise at most $R+2\le3R$ circular edges can be destroyed:
all block interiors survive, and only block boundaries are candidates.
This bounds the loss of original edges as well, even if some are later
recreated. Summing expected losses gives
$\E J(X_t)\ge n-3t r_n(n-2)q$. Apply the variational definition of total
variation to the function $J/n\in[0,1]$ and use
\eqref{eq:adj-uniform}.

For $\Markov_p$, use the number $R$ of switches of the Markov label word,
not the number of hidden adhesive blocks. Here
$R\sim\operatorname{Bin}(n-1,q/2)$. With zero or one switch the increment
is a rotation or the identity. Otherwise at most $R+1$ edges change, and
$(R+1)\mathbf1_{\{R\ge2\}}\le\tfrac32 R(R-1)$. This proves
\eqref{eq:markov-adj-lower}. The coarser estimate
$\E[J(X_t)]\ge n-3tnq$ is also valid for all four families, by counting
broken adhesive bonds.

For the last assertion, preassign independent Bernoulli$(p)$ coins to
each original directed edge in each round, even after that edge is lost.
Whenever the edge is currently an internal packet contact, use its
assigned coin for that bond. This preserves the model's independent
bond law. There are at most two exempt circular boundaries per round
in a cut-first model and at most one in the inverse model. An original
edge whose $t$ coins all survive cannot disappear unless it encounters
an exempt boundary. At most $2t$ distinct original edges can encounter
such a boundary before being lost. Therefore, on this enlarged space,
\[
 J(X_t)\ge Z_t-2t,\qquad Z_t\sim\operatorname{Bin}(n,p^t).
\]
No independence between $Z_t$ and the chosen boundaries is needed.
Under uniform measure, Markov's inequality and
\eqref{eq:adj-uniform} give
$\pi_n(J\ge\ell)\le n/((n-1)\ell)$, proving the claim.
\end{proof}

\begin{corollary}[Uniform order lower bounds]\label{cor:uniform-lower}
As $n\to\infty$, uniformly over $0\le p<1$, all four families satisfy
$t_{\mix}(1/4)=\Omega(L(n,p))$. The inverse Markov family also satisfies
$t_{\mix}(1/4)=\Omega((nq^2)^{-1})$.
\end{corollary}
\begin{proof}
When $p\le1/2$, $L(n,p)=O(\log n)$, so use
\eqref{eq:support-lower}. Suppose $p>1/2$ and write $x=nq$.
If $x$ is bounded, the coarse expected-loss bound in the preceding proof
gives $\Omega(q^{-1})$, as needed. For $x\ge2$, take
$t=\lfloor\log x/(32q)\rfloor$. Since $-\log p\le2q$,
$\mu=\E Z_t\ge n x^{-1/16}$, and uniformly for $x\ge2$,
\[
 \frac{2t}{\mu}\le\frac{\log x}{16x^{15/16}}<\frac14.
\]
Also $\mu\ge n^{15/16}\to\infty$. Binomial concentration, with
$\ell=\lfloor\mu/2\rfloor$, makes the right side of
\eqref{eq:survival-binomial} tend to one. Values of $x$ in a fixed
bounded interval are already covered by the expected-loss bound.
Combining the bounds proves the uniform order assertion. The additional
Markov bound follows directly from \eqref{eq:markov-adj-lower}.
\end{proof}

\subsection{Interpolation retains a stronger small-adhesion obstruction}

Let $T_{n,p}$ be the number of block-key ties in one $\Interp_p$ round,
and define $s_n=r_n(n-2)$. A complete effective seed consists of the cut,
the bond coins actually used, and the fair block-tie coins. Its entropy is
\[
 h_{n,p}=\Ent(K)+s_n H_2(p)+\E T_{n,p}.
\]
Here and below $0\log 0=0$. A useful explicit upper bound is
\begin{equation}\label{eq:sticky-entropy-upper}
 \boxed{\displaystyle
 \bar h_{n,p}=\Ent(K)+s_n H_2(p)
   +\min\!\left\{\E T_n(K)+s_n p,\ r_n+\frac{s_n q}{2}\right\}.}
\end{equation}
The quantities $T_n(k)$ and $\Ent(K)$ are the finite binomial sums
already given for the unglued model.

\begin{proposition}[Parameter-dependent information bound]
\label{prop:sticky-information}
One has $h_{n,p}\le\bar h_{n,p}$. For $0<\varepsilon\le1/2$,
\begin{equation}\label{eq:sticky-fano}
 t_{\mix}^{\Interp_p}(\varepsilon)\ge
 \frac{(1-\varepsilon)\log_2(n!)-H_2(\varepsilon)}{\bar h_{n,p}}.
\end{equation}
For every sequence $p=p_n\in[0,1]$ and fixed $\delta\in(0,1)$,
\begin{equation}\label{eq:sticky-strong}
 d_n\!\left(\left\lfloor(1-\delta)
       \frac{\log_2(n!)}{\bar h_{n,p_n}}\right\rfloor\right)
       \longrightarrow1.
\end{equation}
For even $n$, if $p_n=o((\sqrt n\log n)^{-1})$, the original leading
lower scale $\sqrt{2\pi n}\log_2 n$ is therefore unchanged.
\end{proposition}

\begin{proof}
A tie involving two singleton blocks was already an unglued tie. Every
other tie uses a nonsingleton block, and there are at most as many such
blocks as surviving bonds. Also $T_{n,p}\le B/2$. Taking expectations
proves \eqref{eq:sticky-entropy-upper}. The entropy argument of
Proposition~\ref{prop:finite} then proves \eqref{eq:sticky-fano}.

For completeness, the information-concentration argument remains uniform
in $p$, even when $p$ tends to an endpoint. Write the one-round surprisal
as $Z=Z_K+Z_{\rm bonds}+T_{n,p}$. The cut term has variance $O(n)$ by
Lemma~\ref{lem:information}. A single Bernoulli surprisal has variance
$p q\log_2^2(p/q)$, uniformly bounded for $0\le p\le1$. Conditional on
a nonempty cut there are exactly $n-2$ such variables. Empty cuts add
at most $O(n^2 2^{-n})$ to the variance. Finally
$\Var(T_{n,p})\le(n/2)\E T_{n,p}$. The inequality for the variance of a
sum yields
\[
 \Var(Z)\le Cn(1+\E T_{n,p})\le Cn h_{n,p}.
\]
For $t$ as in \eqref{eq:sticky-strong}, the total seed surprisal has
mean at most $(1-\delta)\log_2(n!)$ and variance at most
$Cn\log_2(n!)$. Chebyshev's inequality makes its upper deviation by
$\delta\log_2(n!)/2$ have probability $O(1/\log n)$.
Histories below this threshold have at most
$(n!)^{1-\delta/2}$ outputs, whose uniform mass tends to zero.
This proves \eqref{eq:sticky-strong}. At the stated small-$p$ scale,
$s_n(H_2(p)+p)=o(\sqrt n)$, so
$\bar h_{n,p}\sim h_n\sim\sqrt{n/(2\pi)}$ for even $n$.
\end{proof}

\section{Local collision bounds: polynomial dependence on deck size}
\label{sec:local-upper}

The proof method is deliberately common to the four families. A local
alternative that exchanges two input cards supplies a random-swap
projection. Three layers of such alternatives connect the deck along a
path. Keeping the full central transposition spectrum then avoids a
factor of order $n$ in the final mixing bound.

\subsection{Extracting fair swaps without requiring a rare global event}

\begin{lemma}[Random-subset projection extraction]\label{lem:extract}
Let $N$ be a matching. Suppose an increment law can be sampled by first
choosing auxiliary data, including $S\subseteq N$, and then independent
fair swaps on the edges of $S$, acting on the input of a fixed remaining
permutation. If every $e\in N$ has $\Prob(e\in S)\ge\rho$, then its
adjoint norm deficit satisfies
\begin{equation}\label{eq:extract}
 \mathcal D(f)\ge\rho\inner f{(I-C_N)f}.
\end{equation}
It suffices that the two local alternatives for an edge have masses
$a,b$: their common part supplies a fair swap with total mass
$2\min(a,b)$; unmatched mass is retained as a separate outcome.
\end{lemma}
\begin{proof}
Conditionally the adjoint is $U_{d^{-1}}C_S$, so the conditional deficit
is $\inner f{(I-C_S)f}$. Mixture variance can only increase it.
The commuting involutions associated to $N$ are simultaneously
diagonalizable. On a nonconstant matching character with nonempty support
$E\subseteq N$, $I-C_S$ equals one precisely when $S\cap E\ne\varnothing$.
This event has probability at least $\rho$, by considering any one
$e\in E$. This proves the quadratic-form inequality. Splitting two
atoms into an equal common mass and residual masses proves the last
assertion. Independence of the activation events is not required.
\end{proof}

In applications with several disjoint local windows, fix the exterior
data and group outcomes into orbits of the window flips. Each flip leaves
the input counts at the beginning and end of its window unchanged.
For card-weighted release the likelihood therefore factors into local
window likelihoods and an unchanged exterior likelihood. For uniform
block interleaving it is constant on each orbit. For Markov labels,
separated windows are conditionally independent given the exterior labels.
These facts justify independent common-part extractions within a layer.
We never demand that all the windows be active simultaneously.

\begin{lemma}[Three path layers]\label{lem:three-layers}
Suppose $N_1,N_2,N_3$ are matchings partitioning the edges of a Hamiltonian
path, and
$\mathcal D(f)\ge\kappa\inner f{(I-C_{N_j})f}$ for each $j$.
Then
\begin{equation}\label{eq:three-layer-result}
 t_{\mix}(\varepsilon)\le
 \left\lceil\frac{2304n(n-1)}{\kappa}\Lam\right\rceil.
\end{equation}
\end{lemma}
\begin{proof}
Let $A_j=C_{N_j}$ and $S=(A_1+A_2+A_3)/3$.
The path comparison of Lemma~\ref{lem:path} still gives
$L_{\RT}\preceq8n\sum_e F_e$, where $F_e=(I-U_{\tau_e})/2$.
On a central eigenspace $V_\lambda$, Lemma~\ref{lem:DL} bounds the
product $A_3A_2A_1$ and gives loss at least
$\lambda\norm f^2/(64n)$. If $a_j=\norm{(I-A_j)f}$, its loss is at most
\[
 \sum_{j=1}^3(a_1+\cdots+a_j)^2
 \le6\sum_{j=1}^3a_j^2=18\inner f{(I-S)f}.
\]
Thus $\mathcal D(f)\ge\kappa\lambda\norm f^2/(1152n)$ on
$V_\lambda$. All relevant operators preserve these spaces.
The spectral sum \eqref{eq:RT-explicit}, at time
$\kappa t/(2304n)$, proves \eqref{eq:three-layer-result}.
\end{proof}

\subsection{The polynomial estimates}

\begin{theorem}[Explicit local upper bounds]\label{thm:local-upper}
For even $n=2m\ge4$ let $\omega_n=2^{-n}\binom n m$. Define
\begin{align}
 \rho_{\Interp}&=q^4\min\{p^2,q^2/2\},\label{eq:rho-I}\\
 \rho_{\Card}&=\frac{q^6}{1024(m+1)},\label{eq:rho-C}\\
 \rho_{\Block}&=\frac{q^6}{32768(1+mpq)\sqrt{m+1}}.\label{eq:rho-B}
\end{align}
For $\mathsf F\in\{\Interp,\Card,\Block\}$ with $\rho_{\mathsf F}>0$,
\begin{equation}\label{eq:local-cut-upper}
 t_{\mix}^{\mathsf F_p}(\varepsilon)\le
 \left\lceil\frac{2304n(n-1)}{\omega_n\rho_{\mathsf F}}\Lam\right\rceil.
\end{equation}
For all $n\ge4$ and $p<1$,
\begin{equation}\label{eq:local-markov-upper}
 t_{\mix}^{\Markov_p}(\varepsilon)\le
 \left\lceil9216n(n-1)q^{-2}\Lam\right\rceil.
\end{equation}
\end{theorem}

\begin{proof}[Common cut-first geometry]
Condition on $K=m$ and denote input cards in the two packets by
$a_1,\ldots,a_m$ and $b_1,\ldots,b_m$. Use the matching
$M=\{\{a_i,b_i\}\}$ and split
$N=\{\{b_i,a_{i+1}\}:1\le i<m\}$ into odd and even $i$.
Together these form the path $a_1,b_1,a_2,b_2,\ldots,a_m,b_m$.
The gadgets for a fixed parity use disjoint two-position windows.
It suffices to prove activation probability at least
$\rho_{\mathsf F}$ for each edge; multiplying by the cut probability
$\omega_n$ and applying Lemma~\ref{lem:three-layers} then gives
\eqref{eq:local-cut-upper}.
\end{proof}

\begin{proof}[Interpolation activations]
If $a_i,b_i$ are singleton blocks, their tied midpoint keys give a fair
swap of $\{a_i,b_i\}$. At most four bonds must break, so the activation
probability is at least $q^4$.

For an edge $\{b_i,a_{i+1}\}$ isolate positions $i,i+1$ in both packets
by breaking their external bonds, an event of probability at least
$q^4$. Write $a=a_i,b=b_i,c=a_{i+1},d=b_{i+1}$.
If the two internal bonds break and both ties put left first, the local
output is $a,b,c,d$, with conditional mass $q^2/4$. If both internal
bonds survive and the tied two-card blocks put left first, it is
$a,c,b,d$, with mass $p^2/2$. These outputs differ by exactly the input
swap $\{b,c\}$. Blocks outside the isolated window have keys outside
its key interval and hence do not interfere. Common-part extraction
activates the swap with probability at least
$2q^4\min(q^2/4,p^2/2)=\rho_{\Interp}$. Disjoint windows have independent
internal choices after the external bonds are fixed. Lemma~\ref{lem:extract}
applies to each layer.
\end{proof}

\begin{proof}[Card-weighted activations]
Conditional on the bond partitions, realize the interleaver by two
independent continuous-time packet processes. A packet with $a$ cards
remaining waits an exponential time of rate $a$, then emits its next
whole block. Competing exponential clocks give exactly the probability
$A/(A+B)$; the construction is a proof device, not a new model parameter.

Isolate a two-position window $i,i+1$ as four singleton cards, requiring
at most six broken bonds. Let $a=m-i+1\ge2$. In one packet let $S$ be
the time at which its prefix of $i-1$ cards has been emitted. The rates
visited before this time are a subset of $a+1,\ldots,m$, whence
\[
 \E S\le\sum_{r=a+1}^m\frac1r\le\log(m/a).
\]
The two prefix times are independent and identically distributed under
the imposed bonds. Divide $[0,2\log(m/a)]$ into bins of width $1/a$,
including an endpoint bin. There are at most $2(m+1)$ bins since
$a\log(m/a)\le m/e$. At least half of the probability lies in this
interval; if its length is zero, both prefix times are zero.
Cauchy--Schwarz gives probability at least $1/[8(m+1)]$ that both times
lie in one bin. The earlier ready packet waits for the later one with
conditional probability at least $e^{-1}$, since its current rate is $a$.
Thus with probability at least $1/[8e(m+1)]$, both window heads are
exposed before either is emitted.

At that moment both packets contain $a$ cards. Each of the local orders
$a_i,b_i,a_{i+1},b_{i+1}$ and $a_i,a_{i+1},b_i,b_{i+1}$ has conditional
probability at least $1/64$. This follows directly by multiplying the
four card-count ratios; the remaining counts differ by at most two, and
the cases $a=2$ only force some of the choices. After these four singleton
releases both orders leave the same state. Their common part therefore
activates $\{b_i,a_{i+1}\}$ with probability at least
$q^6/[256e(m+1)]\ge\rho_{\Card}$.

For $\{a_i,b_i\}$ use a one-position singleton window, at most four
broken bonds, and the two orders $a_i,b_i$ and $b_i,a_i$. The same prefix
argument applies with $a\ge1$, and each order has probability at least
$1/4$. Its activation probability is larger than $\rho_{\Card}$.
For a layer, local flips preserve its good windows and all their entry
and exit counts; the likelihood factorization preceding
Lemma~\ref{lem:three-layers} justifies simultaneous fair extraction.
\end{proof}

\begin{proof}[Block-weighted activations]
Again isolate a two-position window as four singleton cards, at cost
at most $q^6$. Require that the two packets have the same number of
prefix blocks, and also the same number of suffix blocks. Prefix and
suffix partitions are independent. For a binomial variable $Y$ of
variance $v$ and an independent copy $Y'$, concentration in an interval
of at most $5\sqrt{v+1}$ integers gives
\[
 \Prob(Y=Y')=\sum_j\Prob(Y=j)^2\ge\frac1{10\sqrt{v+1}}.
\]
Indeed, Chebyshev at radius $2\sqrt{v+1}$ leaves probability at least
$3/4$, and Cauchy--Schwarz gives the displayed bound. Prefix and suffix
block counts are shifted binomials or deterministic zeros. Their two
agreements therefore have joint probability at least
$1/[100(1+mpq)]$.

On agreement, let each packet have $u$ prefix blocks and $v$ suffix
blocks, and put $b=u+2+v\le m$. Under uniform block interleaving the
chance of reaching the two window heads together is
\[
 \frac{\binom{2u}u\binom{2b-2u}{b-u}}{\binom{2b}b}
 \ge\frac1{4\sqrt{b+1}}\ge\frac1{4\sqrt{m+1}}.
\]
For example, this inequality follows from
$4^j/(2\sqrt{j+1})\le\binom{2j}j\le4^j/\sqrt{j+1}$, including $j=0$.
Each of the two four-card patterns used above has conditional probability
at least $1/64$, now with block counts replacing card counts.
Their probabilities are equal, so their union is already a fair swap.
Multiplying the estimates gives activation probability at least
$q^6/[12800(1+mpq)\sqrt{m+1}]\ge\rho_{\Block}$.
One-position windows similarly give the matching $M$ with at least
this probability. Within a layer, condition on the partitions and the
exterior interleaving. All allowed orientations then have equal weight,
which verifies Lemma~\ref{lem:extract}.
\end{proof}

\begin{proof}[Markov-label activations]
Put $b=q/2$ and $a=1-b$. Partition the neighboring input edges
$\{i,i+1\}$ into three layers according to $i$ modulo three.
Pairs in one layer have at least one intervening fixed label, and are
conditionally independent given the exterior labels. For an interior
pair, require its two exterior labels to agree and its two interior
labels to differ. Exchanging the interior labels $01\leftrightarrow10$
has equal conditional weight and swaps exactly the two adjacent input
cards in the inverse shuffle. The total activation probability is
$2ab^2\ge q^2/4$: the four relevant length-four patterns each have
probability $ab^2/2$.

For a pair touching an endpoint, common-part extraction conditional on
its single exterior neighbor has activation probability $2b^2$, also
at least $q^2/4$. Apply Lemma~\ref{lem:extract} and then
Lemma~\ref{lem:three-layers} with $\kappa=q^2/4$. This proves
\eqref{eq:local-markov-upper}. The adjacent-mark collision mechanism is
related to that used in the stronger fixed-parameter analysis of
\cite{JM15}; the present argument is a direct parameter-explicit
spectral comparison, not an improvement on their $n$ dependence.
\end{proof}

\section{A sharper high-adhesion bound for the cut-first models}
\label{sec:rotation-upper}

The local bounds lose unnecessary powers of $q$ when almost every contact
survives. A different argument uses the frequent rotations and the rare
one-break moves. Unlike the balanced-cut proof, it works for odd decks too.

\begin{theorem}[Rotation and one-break comparison]\label{thm:rotation}
For $n\ge4$ let $\omega_k=2^{-n}\binom nk$ and define
\[
 c_k=\min\{1/3,k/(2n)\},\qquad
 \eta_n=\sum_{k=2}^{n-2}\min\{\omega_k c_k,\omega_{k+1}c_{k+1}\}.
\]
Then $\eta_n>0$ and $\eta_n\to1/4$. For each of
$\Card_p,\Block_p,\Interp_p$, with $0<p<1$, put
\[
 a_{n,p}=\frac{p^{n-2}}{4n},\qquad
 b_{n,p}=\eta_n q p^{n-3}.
\]
An explicit common upper bound is
\begin{equation}\label{eq:rotation-exact}
 t_{\mix}(\varepsilon)\le
 \left\lceil\frac{72n^2(n-1)}{\min\{a_{n,p},b_{n,p}\}}\Lam\right\rceil.
\end{equation}
In particular \eqref{eq:rot-summary} holds with an absolute constant.
\end{theorem}

\begin{proof}
Write $r=(1,2,\ldots,n)$ for a cyclic rotation and
$A=n^{-1}\sum_{j=0}^{n-1}U_{r^j}$ for its invariant-space projection.
The all-surviving-bonds event shows that the increment law at $p$ contains
$p^{n-2}$ times its law at $p=1$. Empty cuts can be split with this same
weight, so the assertion is unconditional.

For $\Block_1$ and $\Interp_1$ the rotation operator is
$(I+\E U_{r^K})/2$, with $K\sim\operatorname{Bin}(n,1/2)$.
For $\Card_1$ it is the same expression with
$K\sim\operatorname{Bin}(n-1,1/2)$: indeed
$\omega_k(n-k)/n=\tfrac12\Prob(\operatorname{Bin}(n-1,1/2)=k)$.
At any nontrivial $n$th root of unity the binomial multiplier has modulus
at most $|\cos(\pi/n)|^{n-1}$. The inequalities
\[
 \cos(\pi/n)\le1-2/n^2,\qquad
 1-e^{-2(n-1)/n^2}\ge\frac1{2n}
\]
imply a squared-norm deficit at least $1/(4n)$ away from $A$.
Mixture variance therefore gives
\begin{equation}\label{eq:rotation-control}
 \mathcal D(f)\ge a_{n,p}\norm{(I-A)f}^2.
\end{equation}

Now let $T_k=(k,k+1,\ldots,n)$, the output-list permutation that moves
card $k$ to the bottom. At cut $k$, use left blocks
$[1,\ldots,k-1],[k]$ and one whole right block. This requires one
broken bond and $n-3$ surviving bonds, of probability $q p^{n-3}$.
The order left-large, right, left-singleton produces $T_k$.
Its conditional probability is one for interpolation, $1/3$ for
block-weighted release, and
$\frac{k}{n}\frac{n-k}{n-k+1}\ge k/(2n)$ for card-weighted release.
Thus every family gives $T_k$ mass at least
$w_k=\omega_k c_k q p^{n-3}$, for $2\le k\le n-1$.

Allocate half of each such atom to each neighboring pair of indices,
discarding unused endpoint allocations. For the pair $T_k,T_{k+1}$,
the conditional variance contributes at least
\[
 \frac{\min\{w_k,w_{k+1}\}}4
 \norm{(I-U_{(k,k+1)})f}^2
 =\min\{w_k,w_{k+1}\}\inner f{F_{(k,k+1)}f},
\]
because $T_kT_{k+1}^{-1}=(k,k+1)$. Summing disjoint allocations gives
\begin{equation}\label{eq:insertion-control}
 \mathcal D(f)\ge b_{n,p}\sum_k v_k\norm{F_{(k,k+1)}f}^2,
 \qquad
 v_k=\frac{\min\{\omega_kc_k,\omega_{k+1}c_{k+1}\}}{\eta_n}.
\end{equation}
Here $v_k\ge0$ and $\sum_kv_k=1$.

On a central eigenspace $V_\lambda$, put $g=Af$ and
$x=\norm{f-g}^2$, $y=\sum_kv_k\norm{F_{(k,k+1)}f}^2$.
Cyclic invariance makes all cyclic neighboring-edge energies of $g$
equal. The path comparison consequently yields
\[
 \lambda\norm g^2\le8n^2\sum_kv_k\norm{F_{(k,k+1)}g}^2.
\]
Since each $F$ is a contraction and $\lambda\le2$, it follows that
$\lambda\norm f^2\le18n^2(x+y)$.
Equations \eqref{eq:rotation-control} and \eqref{eq:insertion-control}
give $\mathcal D(f)\ge\min(a_{n,p},b_{n,p})(x+y)/2$.
The central spectral sum proves \eqref{eq:rotation-exact}.

Finally, on $|k-n/2|\le n^{2/3}$, the adjacent binomial-mass ratios tend
uniformly to one and $c_k\to1/4$. This interval has binomial probability
tending to one. Outside it the total mass tends to zero. Thus
$\eta_n\to1/4$; its finite values for $n\ge4$ are positive.
This proves the absolute-constant order bound.
\end{proof}

\begin{proof}[Proof of Theorem~\ref{thm:parameter-exponents}]
For fixed $n$, \eqref{eq:rotation-exact} is $O_{n,\varepsilon}(q^{-1})$
as $p\uparrow1$, while \eqref{eq:cut-adj-lower} gives the matching
$\Omega_{n,\varepsilon}(q^{-1})$ lower bound. For $\Markov_p$, combine
\eqref{eq:local-markov-upper} with \eqref{eq:markov-adj-lower}.
The numerator $1-\varepsilon-1/(n-1)$ is positive in the stated range.
\end{proof}

\begin{remark}[Why one parameter has two different singular exponents]
In a cut-first model, the imposed cut already creates two packets.
One additional broken contact can produce a move outside the rotation
subgroup. In the inverse-label model, zero or one label switch still
gives only a rotation. At least two switches are necessary. This is a
structural distinction in the models, not evidence for a universal
physical exponent. For $n=2$ the rotation subgroup is all of $\Sn$, so
the quadratic-exponent statement is deliberately not made.
\end{remark}

\section{Small adhesion, regeneration, and a cutoff regime}
\label{sec:small-adhesion}

\begin{proposition}[A comparison that is valid without commutation]
\label{prop:all-singletons}
For even $n$, let $\beta_n=\min\{\omega_{n/2},\omega_{n/2+1}\}$.
For $p<1$,
\[
 t_{\mix}^{\Interp_p}(\varepsilon)\le
 \left\lceil\frac{32768n(n-1)}{\beta_n q^{n-2}}\Lam\right\rceil.
\]
\end{proposition}
\begin{proof}
With probability $q^{n-2}$ all internal bonds at a nonempty cut break.
Splitting empty-cut identity mass by the same factor yields
$Q_p=q^{n-2}Q_0+(1-q^{n-2})R$ for a contraction $R$.
Mixture variance gives
$\mathcal D_p(f)\ge q^{n-2}\mathcal D_0(f)$.
Apply the central comparison of Lemma~\ref{lem:central} and its full
spectral proof. No claim that arbitrary intervening shuffles commute or
can only accelerate mixing is used.
\end{proof}

\begin{proposition}[A GSR regeneration-window bound]\label{prop:refresh}
Let $\mathsf F\in\{\Card,\Block,\Markov\}$, and put
$\nu_{\Card}=\nu_{\Block}=n-2$, $\nu_{\Markov}=n-1$, and
$r=\lceil\log_2(n(n-1))\rceil$. For $n\ge4$ and $p<1$,
\begin{equation}\label{eq:refresh-bound}
 t_{\mix}^{\mathsf F_p}(\varepsilon)\le
 r\left\lceil2q^{-\nu_{\mathsf F}r}\log(1/\varepsilon)\right\rceil.
\end{equation}
\end{proposition}
\begin{proof}
One round contains $q^{\nu_{\mathsf F}}$ times its $p=0$ kernel, which
is GSR or inverse GSR. A run of $r$ such components therefore has weight
$q^{\nu_{\mathsf F}r}$ in the $r$-step kernel.
After $r$ inverse GSR shuffles, independent uniform $r$-bit labels have
all distinct values with probability
$b=(2^r)_n/2^{rn}\ge1-n(n-1)/2^{r+1}\ge1/2$.
Conditional on this event their order is uniform, so the $r$-step GSR
kernel minorizes $b\pi_n$. This standard strong-stationary construction
is described in \cite{AD86,BD92}. Consequently each $r$-step window
has a uniform component of mass at least $q^{\nu_{\mathsf F}r}/2$.
All residual kernels preserve uniform measure. Iterating the
minorization gives \eqref{eq:refresh-bound}. This argument uses whole
windows, not invalid thinning of isolated successful rounds.
\end{proof}

In particular, for any fixed $c>0$, if $p\le c/(n\log n)$ then
$q^{-\nu_{\mathsf F}r}$ stays bounded as $n\to\infty$. The regeneration
bound and the information floor give
$t_{\mix}^{\mathsf F_p}(1/4)=\Theta_c(\log n)$ for all three
GSR-derived families. The next result strengthens this order statement
to the classical cutoff location under a little-$o$ assumption.

\begin{theorem}[Classical cutoff survives sufficiently weak adhesion]
\label{thm:weak-cutoff}
For $\mathsf F\in\{\Card,\Block,\Markov\}$ and all integers $t\ge0$,
\begin{equation}\label{eq:perturbation}
 |d_{n,p}^{\mathsf F}(t)-d_{n,0}^{\mathsf F}(t)|
 \le t(n-1)p.
\end{equation}
In particular, if $p_n=o((n\log n)^{-1})$, these families have the GSR
cutoff at $\tfrac32\log_2 n$, with an $O(1)$ window. On time sequences
where the classical profile converges, their profiles have the same limit.
\end{theorem}
\begin{proof}
Couple each round with its $p=0$ version until the first surviving bond.
At most $n-1$ contacts are tested in a round, and, conditional on no
surviving bond, the interleaving or labeling is precisely the original
GSR rule. A union bound over the tested contacts and rounds gives
\eqref{eq:perturbation}. The error is $o(1)$ throughout any $O(\log n)$
time range under the hypothesis. Transfer the classical total-variation
cutoff and profile from \cite{BD92}; inversion does not change the
Markov family's distance. This does not prove cutoff for fixed $p>0$.
\end{proof}

\begin{table}[ht]
\centering\small
\begin{tabularx}{\textwidth}{@{}lYY@{}}
\toprule
Regime & What is proved & What is not proved\\
\midrule
$p=0$ & GSR cutoff; slow interpolation lower bound & Interpolation's sharp order or cutoff\\[2pt]
$p=O((n\log n)^{-1})$ & Logarithmic order for $\Card,\Block,\Markov$; GSR cutoff under little-$o$ & A fixed-positive-$p$ cutoff theorem\\[2pt]
Fixed $0<p<1$ & Polynomial cut-first bounds; $O_p(\log^4 n)$ for $\Markov$ & Matching orders for the cut-first families\\[2pt]
$n(1-p)=O(1)$ & Improved cut-first upper bound; inverse-model quadratic obstruction & Sharp powers of $n$\\[2pt]
Fixed $n$, $p\uparrow1$ & Exact linear versus quadratic parameter exponents & A large-deck cutoff conclusion\\
\bottomrule
\end{tabularx}
\caption{Different limits answer different questions. A fixed-$n$
adhesion exponent is not a substitute for a uniform-in-$n$ mixing theorem.}
\label{tab:regimes}
\end{table}

\section{Finite calculations and reproducibility}

Table~\ref{tab:entropy} evaluates the exact binomial sums in the entropy
formula. The entropy time $a_n$ is not an estimate proved to equal the
mixing time, especially for small $n$.

\begin{table}[ht]
\centering
\small
\begin{tabular}{rrrrrr}
\toprule
$n$ & $H(K)$ & $\E T_n(K)$ & $h_n$ & $a_n$ & lower bound at $1/4$\\
\midrule
52 & 3.8973 & 3.6462 & 7.5434 & 29.9043 & 23\\
64 & 4.0471 & 4.8910 & 8.9381 & 33.1161 & 25\\
256 & 5.0471 & 8.5872 & 13.6343 & 123.5117 & 93\\
1024 & 6.0471 & 15.4730 & 21.5201 & 407.4797 & 306\\
16384 & 8.0471 & 54.7737 & 62.8208 & 3275.1458 & 2457\\
\bottomrule
\end{tabular}
\caption{Entropy quantities in bits and rounds. The last column is the
ceiling of Proposition~\ref{prop:finite} with $\varepsilon=1/4$, not an
observed mixing time.}
\label{tab:entropy}
\end{table}

The baseline script \texttt{verify\_shuffle.py} performs the following
diagnostics, recorded in \texttt{verification\_results.json}:
\begin{itemize}
\item Exact rational midpoint enumeration checks the tie formula for
8,382 deck-size/cut pairs with $2\le n\le128$.
\item For each of the 63 even sizes from 4 through 128, it checks the
Hamiltonian-cycle identity, the containment of the nearby-cut tie in $M$,
and the noncentral-divisor bound.
\item It enumerates the full transition matrices on $\mathfrak S_4$
and $\mathfrak S_6$, checks double stochasticity, and checks the two
quadratic-form inequalities numerically. The minimum eigenvalues are
nonnegative up to floating-point errors below $10^{-14}$.
\item Direct matrix iteration gives $t_{\mix}(1/4)=3$ for $n=4$ and
$t_{\mix}(1/4)=4$ for $n=6$.
\end{itemize}
These small cases can detect convention or implementation errors. They
provide no empirical evidence for an asymptotic cutoff and are not used
as premises in the proofs.

The combinatorial and entropy checks require only the Python standard
library. The optional full-matrix checks require NumPy. They can be run as
\begin{verbatim}
python3 verify_shuffle.py --max-n 128 --matrix \
    --output verification_results.json
\end{verbatim}
No Monte Carlo estimate, training procedure, or fitted asymptotic exponent
is used in this manuscript.

For the four families, the extension script
\texttt{verify\_\allowbreak adhesive\_\allowbreak shuffles.py}
additionally checks:
\begin{itemize}
\item 140 exact rational increment laws: all four families, $2\le n\le8$,
and $p\in\{0,1/4,1/2,3/4,1\}$. These checks include normalization, the
block versus updated-count descriptions of $\Card_p$, the two descriptions
of $\Markov_p$, the original interpolation law at $p=0$, and endpoint
supports.
\item The conditional binomial block-count law and the interpolation
entropy bound on the same grid; the insertion identity and cyclic Fourier
inequality through $n=128$.
\item All three local-layer quadratic-form inequalities for every family
at $(n,p)=(4,1/4)$, $(4,1/2)$, $(4,9/10)$, and $(6,1/2)$;
and both rotation/insertion inequalities for the three cut-first families
at these points and at $(5,1/2)$ and $(5,9/10)$.
\item Direct full-state iteration for Table~\ref{tab:finite-mixing},
including the distances immediately before and at each reported threshold.
Every displayed comparison includes all four completed arms; no partial
or timed-out arm is silently omitted.
\end{itemize}
Results are recorded in
\texttt{adhesive\_\allowbreak verification\_\allowbreak results.json}.
The transition probabilities are enumerated exactly before conversion to
floating point for matrix operations. Numerical positive-semidefiniteness
checks use tolerance $10^{-10}$ and do not certify the symbolic inequalities.
Run the extended checks with:
\begin{verbatim}
python3 verify_adhesive_shuffles.py --max-n 8 --matrix \
    --output adhesive_verification_results.json
\end{verbatim}

\begin{table}[ht]
\centering
\begin{tabular}{@{}rr rrrr@{}}
\toprule
$n$ & $p$ & $\Card_p$ & $\Block_p$ & $\Interp_p$ & $\Markov_p$\\
\midrule
4 & $0$ & 3 & 3 & 3 & 3\\
4 & $1/2$ & 4 & 4 & 2 & 7\\
4 & $9/10$ & 19 & 18 & 6 & 143\\
\addlinespace
6 & $0$ & 3 & 3 & 4 & 3\\
6 & $1/2$ & 6 & 6 & 3 & 8\\
6 & $9/10$ & 25 & 24 & 10 & 132\\
\bottomrule
\end{tabular}
\caption{Numerically computed $t_{\mix}(1/4)$ in full-shuffle units,
starting from the identity. These are small-deck diagnostics, not
asymptotic evidence. In particular, interpolation need not slow
monotonically as adhesion increases: breaking the perfectly regular
schedule can initially help. This observation is not a physical calibration.}
\label{tab:finite-mixing}
\end{table}

\section{Cutoff, limitations, and directions for improvement}

\subsection{What remains unresolved}

A family has total-variation cutoff at times $s_n$ if, for each fixed
$\delta\in(0,1)$, the distances at $(1-\delta)s_n$ and
$(1+\delta)s_n$ tend respectively to one and zero.
Theorem~\ref{thm:weak-cutoff} proves cutoff for the GSR-derived families
under a vanishing-adhesion hypothesis. It does not extend to a fixed
positive adhesion parameter, nor to the interpolation family.

For \emph{unglued interpolation}, Theorem~\ref{thm:main} proves
the first limit with $s_n=a_n$, but not the second. For every
fixed $\varepsilon\in(0,1)$, they imply
\[
 \liminf_{\substack{n\to\infty\\n\text{ even}}}
 \frac{t_{\mix}^{(n)}(\varepsilon)}{a_n}\ge1,
 \qquad
 t_{\mix}^{(n)}(1/4)=O(n^{5/2}\log n).
\]
Thus the current lower and upper scales differ by a factor of order
$n^2$. We do not know from these arguments whether the true scale is
$\sqrt n\log n$, a larger polynomial scale, or something between them.
Likewise, the fixed-positive-$p$ cut-first bounds do not match in $n$.
The sharp parameter exponents in Theorem~\ref{thm:parameter-exponents}
hold with $n$ fixed, and cannot settle these large-deck questions.

Sufficient input entropy is not sufficient mixing. GSR shuffling is a
useful warning: its input-entropy lower scale is $\log_2 n$, while its
actual total-variation cutoff is $\tfrac32\log_2 n$ \cite{BD92}.
The word ``entropic'' in this manuscript labels a lower-bound clock,
not a proved location of cutoff.

\subsection{Where the upper argument loses information}

The unglued proof deliberately uses only two cuts. Their probabilities are
$\Theta(n^{-1/2})$, and the induced comparison graph is only a cycle.
All other cuts are discarded through a contraction inequality. This
robust simplification is responsible for much of the potential loss.
It also replaces nonreversible deterministic transport by a norm
comparison, which may conceal accelerated spreading over several rounds.

Three concrete improvements suggest themselves.
First, use a window of $\Theta(\sqrt n)$ typical cuts rather than two
individual cuts, and analyze the resulting collection of conjugated
matchings. Second, study a multistep operator in a frame following the
deterministic merges; expansion of these matchings might give stronger
contraction than a one-step comparison. Third, adapt the collision-based
entropy approach of Morris \cite{Morris09} directly to the cut-dependent
matchings. None of these steps is supplied by the present proof.

For the adhesive cut-first models, the local argument uses just the
balanced cut and very conservative probabilities of meeting in selected
singleton windows. The high-adhesion argument instead retains rotations
and one-break insertions. Taking the better estimate reduces parameter
losses, but neither comparison is expected to capture all of the transport
created by typical shuffles. For $\Markov_p$, the literature's
$O_p(\log^4 n)$ result already improves our elementary $n$ dependence
when $p<1$ is fixed; our contribution there is explicit endpoint control.

\subsection{Other deck sizes and other noise scales}

For odd $n$, \emph{unglued} interpolation has no midpoint ties. The rule then has only the cut
randomness, and elementary support counting already gives an order-$n$
lower bound at fixed error. The sharper even-deck analysis is driven by
an event, the exactly balanced cut, that does not exist for odd decks.
No upper theorem for odd-deck $\Interp_0$ is claimed here. In contrast,
Theorem~\ref{thm:rotation} covers all three cut-first models on odd decks
when $0<p<1$, and the Markov-label bound works on either parity.

Continuous perturbations smaller than $n^{-2}$ all define the same
chain. Perturbations large enough to reverse strict cross-packet
comparisons can inject additional randomness. One must then specify
the noise law, its scaling, and whether within-packet order is enforced.
A single common phase per packet is different again: at an equal cut,
its choices are globally correlated rather than independent for each
pair. These variants should not be conflated in a cutoff claim.

\subsection{Conclusion}

One independently resampled contact-survival probability yields several
distinct, mathematically transparent shuffles. Card-weighted selection
must update remaining card counts after every drop; block-weighted
selection, midpoint scheduling, and correlated-label inversion are
different rules, not equivalent implementations.

The paper establishes explicit parameter-dependent upper and lower
bounds, an interpolation information obstruction, and exact fixed-deck
adhesion exponents: linear for cut-first release and quadratic for
correlated labels. Classical cutoff transfers in a sufficiently weak
perturbation regime. Sharp large-deck mixing orders and cutoff for the
new fixed-parameter cut-first families remain open within this analysis.
The contact model is physically motivated, but is neither a force law
nor an experimentally validated model of hand shuffling.

\end{document}